\documentclass[10pt]{amsart}

\usepackage[all]{xypic}

\usepackage{latexsym}

\usepackage{amssymb}
\usepackage{amsfonts}
\usepackage{amscd}
\usepackage{amsmath,amsthm}

\newtheorem{lemma}{Lemma}[section]

\newtheorem{lem}[lemma]{Lemma}
\newtheorem{prop}[lemma]{Proposition}
\newtheorem{thm}[lemma]{Theorem}
\newtheorem{cor}[lemma]{Corollary}

{

}

\theoremstyle{definition}

\theoremstyle{remark}

\def\BB{{\mathbb B}}

\numberwithin{equation}{section}

\newenvironment{pf}{\noindent{\bf Proof.}}{\hfill $\square$\medskip}

\def\BB{{\mathbb B}}
\def\CC{{\mathbb C}}

\def\PP{{\mathbb P}}

\def\ZZ{{\mathbb Z}}

\def\0ol{{\bar 0}}
\def\1ol{{\bar 1}}
\def\2ol{{\bar 2}}
\def\ol2{{\bar 2}}
\def\3ol{{\bar 3}}
\def\4ol{{\bar 4}}
\def\5ol{{\bar 5}}
\def\6ol{{\bar 6}}
\def\7ol{{\bar 7}}
\def\8ol{{\bar 8}}
\def\9ol{{\bar 9}}

\def\bold0{{\bf 0}}
\def\bold1{{\bf 1}}
\def\bold2{{\bf 2}} 
\def\bold3{{\bf  3}}
\def\bold4{{\bf 4}}
\def\bold5{{\bf 5}}
\def\bold6{{\bf 6}}
\def\bold7{{\bf 7}}
\def\bold8{{\bf 8}}
\def\bold9{{\bf 9}}

\def\P2Skly{\PP^2_{Skly}}

\def\th{\operatorname {th}}    

\def\Aut{\operatorname{Aut}}

\def\dim{\operatorname{dim}}

\def\Fdim{{\sf Fdim}}

\def\Gr{{\sf Gr}}

\def\id{\operatorname{id}}

\def\Pic{\operatorname{Pic}}

\def\Spec{\operatorname{Spec}}

\def\ul1{\operatorname{\underline{1}}}

\def\G{\mathop{\underline{\underline{\it \Gamma}}}\nolimits}

\def\d{\downarrow}

\def\d{\delta}

\def\s{\sigma}

\def\G{\Gamma}

\def\sT{{\sf T}}

\def\cal{\mathcal}

\def\cL{{\cal L}}
\def\cM{{\cal M}}
\def\cN{{\cal N}}
\def\cO{{\cal O}}

\def\Qcoh{{\sf Qcoh}}

\def\dirlim{\mathop{\vtop{\baselineskip -100pt\lineskip -1pt\lineskiplimit 0pt
\setbox0\hbox{lim}\copy0\hbox to \wd0{\rightarrowfill}}}\limits}
\def\invlim{\mathop{\vtop{\baselineskip -100pt\lineskip -1pt\lineskiplimit 0pt
\setbox0\hbox{lim}\copy0\hbox to \wd0{\leftarrowfill}}}\limits}

\def\I11{{1 \kern -0.8pt \! \mbox{l}}} 
\def\mumu{{\mu\kern-4.2pt\mu}}
\def\bfmu{{\mu\kern-4.2pt\mu}}
\def\2slash{\backslash \! \backslash}

\def\boxtimes{\setbox0\hbox{$\Box$}\copy0\kern-\wd0\hbox{$\times$}}

\begin{document}

\title[A NC-HCR for the degree six del Pezzo surface]{A non-commutative homogeneous coordinate ring for the degree six del Pezzo surface}

\author{S. Paul Smith}
\address{ Department of Mathematics, Box 354350, Univ.  Washington, Seattle, WA 98195}
\email{smith@math.washington.edu}



\thanks{The author was supported by the 
National Science Foundation, Award No. 0602347}

\begin{abstract}
Let $R$  be the free $\CC$-algebra on $x$ and $y$ modulo the relations $x^5=yxy$ and $y^2=xyx$
endowed with the $\ZZ$-grading $\deg x=1$ and $\deg y=2$. The ring $R$ appears, in somewhat hidden guise, in a paper on quiver gauge theories.
Let $\BB_3$ denote the blow up of $\CC\PP^2$ at three non-colinear points. The main result in this paper
is that the category of 
quasi-coherent $\cO_{\BB_3}$-modules is equivalent to the quotient of the category of $\ZZ$-graded $R$-modules modulo the full subcategory of modules that are the sum of their finite dimensional submodules. This reduces almost all representation-theoretic questions
about $R$ to algebraic geometric questions about the del Pezzo surface $\BB_3$. For example, the generic simple $R$-module has dimension six. Furthermore,  the main result combined with results of  Artin, Tate, and Van den Bergh, implies  that $R$ is a noetherian domain of global dimension three.
\end{abstract}

\maketitle
 
\section{Introduction}

We will work over the field of complex numbers. 

\subsection{}
The surface obtained by blowing up $\PP^2$ at three non-colinear points is, up to isomorphism,
 independent of the points. It is called the {\sf del Pezzo surface of degree six} and we will denote it by $\BB_3$. 
 
 \subsection{}
Let $R$ be the free $\CC$-algebra on $x$ and $y$ modulo the relations  
\begin{equation}
\label{relations}
x^5=yxy \qquad \hbox{and } \qquad y^2=xyx.
\end{equation}
Give $R$ a $\ZZ$-grading by declaring that 
$$
\deg x=1 \qquad \hbox{and } \qquad \deg y=2.
$$

The ring $R$ arises, in somewhat hidden guise, in a paper about string theory \cite{BP} (see section 
\ref{sect.mtvn}). The present paper concerns only the mathematical properties of $R$ and its relation to the 
degree 6 del Pezzo surface.

\subsection{}
The main result in this paper establishes the following surprising relationship between $R$ and the degree six del Pezzo surface.  
 
 \begin{thm}
 \label{thm.main1}
 Let $R$ be the non-commutative algebra $\CC[x,y]$ defined by the relations (\ref{relations}). 
 Let $\Gr R$ be the category of $\ZZ$-graded left $R$-modules. There is an equivalence of categories 
$$
\Qcoh \BB_3 \equiv {{\Gr R}\over{\Fdim R}}
$$
where the left-hand side is the category of quasi-coherent $\cO_{\BB_3}$-modules 
and the right-hand side is the quotient   category modulo
the full subcategory $\Fdim R$ consisting of those modules that are the sum of their finite dimensional submodules.
\end{thm}
 
 Theorem \ref{thm.main1} is a consequence of the following result.
 
 \begin{thm}
\label{thm.main2}
Let $R$ be the non-commutative algebra $\CC[x,y]$ defined by the relations (\ref{relations}). 
Let $\cL=\cO(-E)$ be the invertible $\cO_{\BB_3}$-module corresponding to a 
$(-1)$-curve $E$ and $\s$ an order 6 
automorphism of $\BB_3$ that cyclically permutes the six  $(-1)$-curves on $\BB_3$.
Then $R$ is isomorphic to the twisted homogeneous coordinate ring
$$
B(\BB_3,\cL,\s) : = \bigoplus_{n \ge 0} H^0(\BB_3,\cL_n)
$$
where
$$
\cL_n:= \cL \otimes (\s^*)\cL \otimes \cdots \otimes (\s^*)^{n-1}\cL.
$$
\end{thm}

In the terminology of Artin, Tate, and Van den Bergh \cite{ATV1} and Artin and Van den Bergh \cite{AV}, $B(\BB_3,\cL,\s)$ is a {\it twisted homogeneous coordinate ring} of $\BB_3$.
Results of Artin, Tate, and Van den Bergh, and Stephenson \cite{St1} now imply that $R$ is a 
3-dimensional Artin-Schelter regular algebra and therefore has the following properties.

\begin{cor}
Let $R$ be the non-commutative algebra $\CC[x,y]$ defined by the relations (\ref{relations}). Then
\begin{enumerate}
  \item 
  $R$ is a left and right noetherian domain;
  \item 
  $R$ has global homological dimension 3;
  \item{}
  $R$ is Auslander-Gorenstein and Cohen-Macaulay in the non-commutative sense;
  \item 
the Hilbert series of $R$ is the same as that of the weighted polynomial ring on three variables 
of weights 1, 2, and 3;
  \item{}
  $R$ is a finitely generated module over its center \cite[Cor. 2.3]{ST94};
  \item{}
  $R^{(6)}:=\oplus_{n=0}^\infty R_{6n}$ is isomorphic to $\bigoplus_{n=0}^\infty H^0(\BB_3,\cO(-nK)$ where $K=K_{\BB_3}$ is the canonical divisor on $\BB_3$;
  \item{}
  $\Spec R^{(6)}$ is the canonical cone over $\BB_3$, i.e., the cone obtained by collapsing the zero section of the total space of the canonical bundle over $\BB_3$.
\end{enumerate}
\end{cor}

This close connection between $R$ and $\BB_3$ means that almost all aspects of the representation theory
of $R$ can be expressed in terms of the geometry of $\BB_3$. We plan to address  this question in another 
paper.

\subsection{}
The justification for calling $R$ a non-commutative homogeneous coordinate ring for $\BB_3$ is the similarity between the equivalence of categories in Theorem \ref{thm.main1} and following theorem of Serre:
\begin{quote}
 if $X \subset \PP^n$ is the scheme-theoretic   
zero locus of a graded ideal $I$ in the polynomial ring $S=\CC[x_0,\ldots,x_n]$ with its standard grading,
and $A=S/I$, then there is an equivalence of categories
\begin{equation}
\label{eq.Serre}
\Qcoh X \equiv {{\Gr A}\over{\Fdim A}}
\end{equation}
where the right-hand side is the quotient category of $\Gr A$, the category of graded $A$-modules,
by the full subcategory $\Fdim A$ consisting of modules whose non-zero finitely generated submodules have
support only at the origin. 
\end{quote}

 \subsection{Motivation}
 \label{sect.mtvn}
 The results in this paper are a prerequisite for some results in \cite{Sm}
 where three superpotential algebras appearing in the string theory literature are investigated by
 relating them to twisted homogeneous coordinate rings In \cite {BP}, Beasley and Plesser study a superpotential algebra they dub the $dP_3I$ path algebra. In 
 \cite{Sm}, we will show that the $dP_3I$ path algebra is isomorphic to $R \rtimes \mu_6$, the 
 skew group ring for  the $6^{\th}$ roots of unity acting on $R$ by $\xi \cdot r=\xi^n r$ for $r \in R_n$; the isomorphism is established in \cite{Sm}. An intimate understanding of $R$ therefore leads to an detailed understanding of the $dP_3I$ path algebra. The $dP_3$ in the notation $dP_3I$ refers to the de Pezzo surface 
 obtained by blowing up 3 non-colinear points in $\PP^2$. The I  in $dP_3I$ is to distinguish this algebra from  two other path algebras with relations that Beasley and Plesser associate to the degree-six del Pezzo surface.

\section{$R=\CC[x,y]$ with $x^5=yxy$ and $y^2=xyx$ is an iterated Ore extension}
\label{sect.nc.ring}

The following result is a straightforward calculation. The main point of it is to show that $R$ 
has the same Hilbert series as the weighted polynomial ring on three variables of weights 1, 2, and 3.

\begin{prop}
[Stephenson \cite{St1},\cite{St2}]
\label{prop.steph}
The ring $R:=\CC[x,y]$ with defining relations 
$$
x^5=yxy \qquad \hbox{and} \qquad  y^2=xyx
$$
is an iterated Ore extension of the polynomial ring $\CC[w]$. Explicitly, if $\zeta$ is a fixed primitive $6^{\th}$ root of unity,  $R$ has the following properties.
\begin{enumerate}
  \item 
   $R= \CC[w][z;\s][x;\tau,\d]$  where $\s \in \Aut \CC[z]$, $\tau \in \Aut \CC[w][z;\s]$, and $\d$ is a $\tau$-derivation defined as follows: 
   \begin{align*}
   \s(w)& =\zeta w,
   \\
   \tau(w) &=  - \zeta^2 w, \qquad    \tau(z)  = \zeta z,
   \\
   \d(w) & = z, \quad  \phantom{xxxxx|}  \d(z)  =- w^2.
   \end{align*}
  \item 
  A set of defining relations of $R=\CC[z,w,x]$ is given by 
  \begin{align*}
zw = & \zeta wz,
\\
xw = & -\zeta^2wx+z ,
\\
xz  = & \zeta zx- w^2.
\end{align*}
  \item 
  $R$ has basis $\{ w^iz^jx^k \; | \; i,j,k \ge 0\}$.
  \item{}
  $R$ is a noetherian domain.
  \item{}
  The Hilbert series of $R$ is  $(1-t)^{-1}(1-t^2)^{-1}(1-t^3)^{-1}$.
\end{enumerate}
\end{prop}
\begin{pf}
Define the elements
\begin{align*}
w: & = \,  y-x^2
\\
z: & = \, xw + \zeta^2 wx 
\\
&= \, xy +\zeta^2 yx - \zeta x^3
\end{align*}
of $R$. Since $y$ belongs to the subalgebra of $R$ generated by $x$ and $w$, $\CC[x,y]=\CC[x,w]=\CC[x,w,z]$. 
It is easy to check that 
\begin{equation}
\label{steph.relns}
zw=\zeta wz, \quad xw= z- \zeta^2 wx, \quad xz=\zeta zx - w^2.
\end{equation}

Let $R'$ be the free algebra $\CC\langle w,x,z \rangle$ modulo the relations
 in (\ref{steph.relns}). We will show $R'$ is isomorphic to $R$. We already know there is a 
 homomorphism $R' \to R$ and we will now exhibit a homomorphism $R \to R'$ by showing there are elements $x$ and $Y$ in $R'$ that satisfy the defining relations for $R$.
 Define the element $Y:=w+x^2$ in $R'$.
A straightforward computation in $R'$ gives 
 $$
 xwx-x^2w = w^2+wx^2
 $$
 so
 $$
 Y^2= w^2+x^2w + wx^2 + x^4 =  xwx+x^4=xYx.
 $$
The next calculation uses the identity $1-\xi+\xi^2=0$ repeatedly. Deep breath...
 \begin{align*}
YxY& = (w+x^2)xw  +wx^3 + x^5
 \\
 & = (w+x^2) (z-\zeta^2wx)+ \big[ wx^3 + x^5 \big]
\\
& = x^2 z -\zeta^2 x^2wx + \big[  wz -\zeta^2 w^2x +wx^3 + x^5 \big]
\\
& = x(\zeta zx -w^2)  -\zeta^2 x(z-\zeta^2wx)x +  \big[  wz -\zeta^2 w^2x +wx^3 + x^5 \big]
\\
& =(\zeta-\zeta^2)xzx - xw^2 - \zeta xwx^2 +  \big[  wz -\zeta^2 w^2x +wx^3 + x^5 \big]
\\
& =(\zeta-\zeta^2)(\zeta zx -w^2)x  -   (z-\zeta^2wx) w    -\zeta(z- \zeta^2 wx)x^2 
\\
& \qquad  \qquad + \big[  wz -\zeta^2 w^2x +wx^3 + x^5 \big]
\\
& =(\zeta^2-\zeta^3)zx^2 -(\zeta-\zeta^2)w^2x - zw +\zeta^2 wxw -\zeta zx^2 - wx^3  
\\
& \qquad  \qquad + \big[  wz -\zeta^2 w^2x +wx^3 + x^5 \big]
\\
& =(\zeta^2-\zeta^3-\zeta)zx^2   +\zeta^2 wxw  + \big[ (1-\zeta) wz -\zeta w^2x + x^5 \big]
\\
& = \zeta^2 wxw   +   \big[ (1-\zeta) wz -\zeta w^2x + x^5 \big]
 \\
& = \zeta^2 w(z-\zeta^2wx)   + \big[ -\zeta^2  wz -\zeta w^2x + x^5 \big]
\\
& = x^5.
 \end{align*}
Since $YxY=x^5$, $R$ is isomorphic to $R'$. Hence $R$ is an
 iterated  Ore extension as claimed. The other parts of the proposition follow easily.   
\end{pf}

It is an immediate consequence of the relations that $x^6=y^3$. Hence $x^6$ is in the center of $R$.

 \section{The del Pezzo surface $\BB_3$}

 Let $\BB_3$ be the surface obtained by blowing up the complex projective plane $\PP^2$ at three 
 non-collinear points. We will write 
 $$
 \pi:\BB_3 \to \PP^2
 $$
 for the morphism that contracts the exceptional curves $E_1$, $E_2$, and $E_3$.
The $(-1)$-curves on $\BB_3$ lie in the following configuration
 \begin{equation}
 \label{six.lines}
 \qquad
\UseComputerModernTips
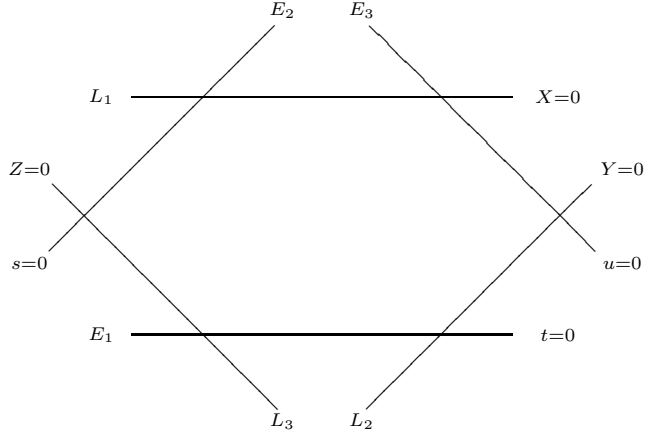
\xymatrix{
&& &  \save []+<0cm,0.1cm>*\txt<4pc>{$\scriptstyle{E_2}$} \restore \ar@{-}[dddlll] 
& \ar@{-}[dddrrr]   \save []+<0cm,0.1cm>*\txt<4pc>{$\scriptstyle{E_3}$} \restore 
\\ 
& \save []+<-0.3cm,0cm>*\txt<4pc>{$\scriptstyle{L_1}$} \restore 
\ar@{-}[rrrrr] &&&&& \save []+<0.5cm,0cm>*\txt<4pc>{$\scriptstyle{X=0}$} \restore  &
\\  
\save []+<-0.2cm,0.1cm>*\txt<4pc>{$\scriptstyle{Z=0}$} \restore  &&&&&&&
\save []+<0.3cm,0.1cm>*\txt<4pc>{$\scriptstyle{Y=0}$} \restore 
&&&&& 
\\ 
\save []+<-0.2cm,-0.1cm>*\txt<4pc>{$\scriptstyle{s=0}$} \restore   &&&&&&&  \save []+<0.3cm,-0.1cm>*\txt<4pc>{$\scriptstyle{u=0}$} \restore    &
\\
&
 \save []+<-0.3cm,0cm>*\txt<4pc>{$\scriptstyle{E_1}$} \restore 
 \ar@{-}[rrrrr] &&&&&  \save []+<0.5cm,0cm>*\txt<4pc>{$\scriptstyle{t=0}$} \restore 
\\
&&& \ar@{-}[uuulll]   \save []+<0cm,-0.1cm>*\txt<4pc>{$\scriptstyle{L_3}$} \restore  
&  \save []+<0cm,-0.1cm>*\txt<4pc>{$\scriptstyle{L_2}$} \restore  \ar@{-}[uuurrr]  }
\end{equation}
where $L_1$, $L_2$, and $L_3$ are the strict transforms of the lines in $\PP^2$ spanned by the points that are blown up. (The labeling of the equations for the $(-1)$-curves is explained in section \ref{sect.curves}.)

The union of the $(-1)$-curves is an anti-canonical divisor so we write
$$
-K:=L_1+L_2+L_3+E_1+E_2+E_3 
$$
($K$ for kanonical). This is, of course, an ample divisor.

 \subsection{The Picard group of $\BB_3$}
 \label{sect.Pic}
 The morphism $\pi:\BB_3 \to \PP^2$ induces an injective group homomorphism
 $\pi^* : \Pic \PP^2 \to \Pic \BB_3$. We write $H=\pi^*L$ where $L$ is a line in $\PP^2$. Hence
 $$
 \Pic \BB_3 = \ZZ H \oplus \ZZ E_1 \oplus \ZZ E_2 \oplus \ZZ E_3.
 $$
 We identify $ \Pic \BB_3$ with $\ZZ^4$ by using the ordered basis
 $$
 H,\; -E_2, \; -E_1, \; -E_3.
 $$
Thus 
$$
H=(1,0,0,0), \quad E_1=(0,0,-1,0), \quad E_2=(0,-1,0,0), \quad E_3=(0,0,0,-1).
$$
In this basis the anti-canonical divisor is
$$
-K=(3,1,1,1).
$$
The Picard group may be presented more symmetrically as
 $$
 \Pic \BB_3 = {{\oplus_{i=1}^3 (\ZZ L_i \oplus \ZZ E_i)}\over{(E_i+L_j=E_j+L_i \; | \; 1\le i,j \le 3)}}.
 $$
 It follows that 
 $$
 H=L_1+E_2+E_3=L_2+E_1+E_3=L_3+E_1+E_2
 $$
 and 
 $$
 L_1=(1,1,0,1), \quad L_2=(1,0,1,1), \quad L_3=(1,1,1,0).
 $$

 \subsection{Cox's homogeneous coordinate ring}
 
 By definition, Cox's homogeneous coordinate ring \cite{Co} for a complete smooth toric variety $X$ is
 $$
S:= \bigoplus_{[\cL] \in \Pic X} H^0(X,\cL).
 $$
From now on, $S$ denotes Cox's homogeneous coordinate ring for $\BB_3$.
  
Let $X,Y,Z,s,t,u$ be coordinate functions on $\CC^6$. 
One can present $\BB_3$ as a toric variety by defining it as the orbit space
$$
\BB_3:={{\CC^6 -W}\over{(\CC^\times)^4}}
$$
where the irrelevant locus, $W$, is the union of nine codimension two subspaces, namely
\begin{equation}
\label{irrel.locus}
\begin{array}{ccc}
\begin{array}{c}
X=t=0\\
Y=s=0\\
Z=u=0
\end{array}
\qquad
\begin{array}{c}
X=Y=0\\
Y=Z=0\\
Z=X=0
\end{array}
\qquad
\begin{array}{c}
s=t=0\\
u=t=0\\
s=u=0
\end{array}
\end{array}
\end{equation}
 and $(\CC^\times)^4$ acts with weights 
$$
\begin{array}{rrrrrr}
X & \phantom{x} & 1 & 1 & 0 & 1 \\
Y &  \phantom{x} & 1 & 0 & 1 & 1 \\
Z & \phantom{x} & 1 & 1 & 1& 0 \\
s & \phantom{x} & 0 & -1 & 0 & 0\\
t & \phantom{x} & 0 & 0 & -1 & 0 \\
u & \phantom{x} & 0 & 0 & 0 & -1
\end{array}
$$
Therefore $S$ is the $\ZZ^4$-graded polynomial ring 
 $$
 S=\CC[X,Y,Z,s,t,u]
 $$
 with the  degrees of the generators given by their weights under the $(\CC^\times)^4$ action, 
 e.g., $\deg X=(1,1,0,1)$, $\deg u=(0,0,0,-1)$, etc.
 It follows from Cox's  results \cite[Sect. 3]{Co} that
$$
\Qcoh \BB_3 \equiv {{\Gr(S,\ZZ^4)}\over{\sT}}
$$
where $\Gr(S,\ZZ^4)$ is the category of $\ZZ^4$-graded $S$-modules and $\sT$ is the full subcategory 
consisting of the modules that are sums of finitely generated $S$-modules supported on $W$.

\subsection{}
\label{sect.curves}
The labelling of the $(-1)$-curves in the diagram (\ref{six.lines}) is explained by the existence of the morphisms
$$
\UseComputerModernTips
\xymatrix{
& \BB_3 \ar[dl]_{\pi} \ar[dr]  
\\
\PP^2 && \PP^2 
}
\qquad
\UseComputerModernTips
\xymatrix{
&  \save []+<-0.2cm,0.2cm>*\txt<4pc>{$\text{ $(X,Y,Z,s,t,u) $}$} \restore     \ar[dl]_{\pi} \ar[dr] 
\\
 \save []+<-0.6cm,-0.2cm>*\txt<4pc>{$\text{$(Xsu,Ytu,Zst)$}$} \restore&&
 \save []+<+0.3cm,-0.2cm>*\txt<4pc>{$\text{$(YZt,XZs,XYu)$}$} \restore
}
$$
that collapse the  $(-1)$-curves: for example,
$\pi$ contracts the three divisors $t=0$, $s=0$, and $u=0$, i.e., $E_1$, $E_2$, and $E_3$.

\subsection{An order six automorphism $\s$ of $\BB_3$}

The cyclic permutation of the six $(-1)$-curves on $\BB_3$ extends to a global automorphism of $\BB_3$
of order six. We now make this explicit.

The category of graded rings consists of pairs $(A,\G)$ consisting of an abelian group $\G$ and a $\G$-graded ring $A$. A morphism $(f,\theta):(A,\G) \to (B,\Upsilon)$ consists of a ring homomorphism $f:A \to B$
and a group homomorphism $\theta:\G \to \Upsilon$ such that $f(A_i) \subset B_{\theta(i)}$ for all $i \in \G$.

Let  $\tau:S \to S$ be the automorphism induced by the cyclic permutation
\begin{equation}
\label{defn.tau}
\UseComputerModernTips
\xymatrix{
X \ar[r] & u \ar[r] & Y \ar[r]_\tau & t \ar[r] & Z \ar[r] & s \ar@{}@/_1pc/[lllll] | {<}
}
\end{equation}
and let $\theta:\ZZ^4 \to \ZZ^4$ be left multiplication by the matrix
$$
\theta= \begin{pmatrix}
2 & -1 & -1 & -1   \\
1 & -1 & -1 & 0 \\
1 & 0 & -1 & -1 \\
1 & -1 & 0 & -1
\end{pmatrix}.
$$
Then $(\tau,\theta):(S,\ZZ^4) \to  (S,\ZZ^4)$ is an automorphism in the category of graded rings.
The irrelevant locus (\ref{irrel.locus}) is stable under the action of $\tau$, so $\tau$ induces an automorphism $\s$ of $\BB_3$.  It  follows from the definition of $\tau$ in (\ref{defn.tau})
 that $\s$ cyclically permutes the six $-1$-curves.

Since $(\tau,\theta)^6=\id_{(S,\ZZ^4)}$ the order of $\s$ divides six. But the action of 
$\s$ on the set of $(-1)$-curves has order six, so $\s$ has order six as an automorphism of $\BB_3$.

\subsection{}

Fix a primitive cube root of unity $\omega$. 
The left action of $\theta$ on $\ZZ^4=\Pic \BB^3$ has  eigenvectors
$$
v_1=\begin{pmatrix} 1    \\  1 \\ 1 \\ 1 \\ \end{pmatrix}, \quad
v_2=\begin{pmatrix} 3    \\  1 \\ 1 \\ 1 \\ \end{pmatrix}, \quad
v_3=\begin{pmatrix} 0    \\  1 \\ \omega \\ \omega^2 \\ \end{pmatrix}, \quad
v_4=\begin{pmatrix} 0    \\  1 \\ \omega^2 \\ \omega \\ \end{pmatrix},
$$
with corresponding eigenvalues $-1$, $+1$, $\omega^2$, $\omega$.

\section{A twisted hcr for $\BB_3$}

In this section we will prove the main theorem: $R$ is isomorphic to a twisted homogeneous coordinate ring
$B=B(\BB_3,\cL,\s)$ for $\BB_3$. The degree-$n$ homogeneous component of $B$ is the global sections of an invertible $\cO_{\BB_3}$-module $\cL_n$; i.e., $B_n=H^0(\BB_3,\cL_n)$.
After defining $\cL_n$  in section \ref{sect.Dn} we prove some vanishing results for
its cohomology that will be used later  to prove that $B$ is generated as a $\CC$-algebra by $B_1$ and $B_2$. Since $R$ is generated as a $\CC$-algebra by $x\in R_1$ and $y \in R_2$ this will allow us to prove that the homomorphism $\Phi:R \to B$ defined in Proposition \ref{prop.R->B} is surjective. We also compute
$\dim H^0(\BB_3,\cL_n) = \dim B_n$ and observe that this is the same as $\dim R_n$ which allows us to conclude that $\Phi$ is an isomorphism.

We write $K$ for the canonical divisor on $\BB_3$.
 
\subsection{A sequence of line bundles on $\BB_3$}
 \label{sect.Dn}
 
We will blur the distinction between a divisor $D$ and the class of the line
bundle $\cO(D)$  in $\Pic\BB_3$.

We define a sequence of divisors:  $D_0$ is zero; $D_1$ is the line $L_1$; for $n \ge 1$ 
$$
D_n:=(1+\theta + \cdots + \theta^{n-1})(D_1).
$$
We will write $\cL_n:=\cO(D_n)$. Therefore 
$$
 \cL_n  = \cL_1 \otimes \s^* \cL_1 \otimes \ldots \otimes (\s^*)^{n-1} \cL_1.
$$
For example,
$$
\begin{array}{cc}
\begin{array}{l}
\cO(D_1)=\cL_1 = \cO(1,1,0,1) \\
\cO(D_2)=\cL_2 = \cO(1,1,0,0) \\
\cO(D_3)=\cL_3 = \cO(2,1,1,1) \\
\cO(D_4)=\cL_4 = \cO(2,1,0,1) \\
\cO(D_5)=\cL_5 =  \cO(3,2,1,1) \\
\cO(D_6)=\cL_6 = \cO(3,1,1,1) \\
\phantom{xyz} \\
\end{array}
\begin{array}{l} 
=\cO(L_1 )\\
=\cO(L_1 + E_3 )\\
=\cO(L_1 + E_3+ L_2  )  \\
=\cO(L_1 + E_3+ L_2 + E_1)   \\
=\cO(L_1 + E_3+ L_2 + E_1 + L_3)  \\
=\cO(L_1 + E_3+ L_2 + E_1 + L_3 +E_2) \\
=\cO(-K).
\end{array}
\end{array}
$$

\begin{lem}
\label{lem.Dn}
Suppose $m \ge 0$ and $0 \le r \le 5$. Then
$$
D_{6m+r} = D_{r}-mK.
$$
\end{lem}
\begin{pf}
Since $\theta^6=1$, 
\begin{align*}
\sum_{i=0}^{6m+r-1}\theta^i & =(1+\theta+\cdots+\theta^5)\sum_{j=0}^{m-1}\theta^{6j} + \theta^{6m}(1+\theta + \cdots+\theta^{r-1})
\\
&=(1+\theta + \cdots+\theta^{r-1})  +  m(1+\theta+\cdots+\theta^5)
\end{align*}
where the sum $1+\theta+\cdots +\theta^{r-1}$ is empty and therefore equal to zero when $r=0$. 
Therefore  $D_{6m+r}= D_{r}  + mD_6   = D_{r}-mK$, as claimed. 
\end{pf}

\subsection{Vanishing results}

For a divisor $D$ on a smooth surface $X$,   we write $$h^i(D):=\dim H^i(X,\cO_{X}(D)).$$
We need to know that $h^1(D)=h^2(D)=0$ for various divisors $D$ on $\BB_3$.  

If $D-K$ is ample, then the Kodaira Vanishing Theorem implies that $h^0(K-D)=h^1(K-D)=0$ and Serre duality then gives $h^2(D)=h^1(D)=0$.

The notational conventions in section \ref{sect.Pic} identify $\Pic \BB_3$ with $\ZZ^4$ via
$$
aH-cE_1-bE_2-dE_3 \equiv (a,b,c,d).
$$
 The intersection form on $\BB_3$ is given by
$$
H^2=1, \qquad E_i.E_j=-\d_{ij}, \qquad H.E_i=0,
$$
so the induced intersection form on $\ZZ^4$ is
$$
(a,b,c,d) \cdot(a',b',c',d')=aa'-bb'-cc'-dd'.
$$

 \begin{lem}
 \label{lem.van}
Let $D=(a,b,c,d) \in \Pic\BB_3\equiv \ZZ^4$.  Suppose that 
\begin{equation}
\label{eq.ample2}
(a+3)^2> (b+1)^2+(c+1)^2+(d+1)^2
\end{equation}
and 
\begin{equation}
\label{eq.ample3}
b,\, c,\, d >-1, \qquad \hbox{and} \qquad  a+1>b+c,\, b+d, \, c+d.
\end{equation}
Then $D-K$ is ample, whence $h^1(D)=h^2(D)=0$.
\end{lem}
\begin{pf}
The effective cone is generated by $L_1$, $L_2$, $L_3$, $E_1$, $E_2$, and $E_3$  so, by the  Nakai-Moishezon criterion, 
$D-K$ is ample if and only if $(D-K)^2>0$ and $(D-K).L_i>0$ and $(D-K).E_i>0$ for all $i$. Now 
$D-K=(a+3,b+1,c+1,d+1)$, so $(D-K)^2>0$ if and only if (\ref{eq.ample2}) holds and $(D-K).D'>0$
for all effective $D'$ if and only if (\ref{eq.ample3}) holds. 

Hence the hypothesis that  (\ref{eq.ample2}) and  (\ref{eq.ample3}) hold implies that $D-K$ is ample. 
The Kodaira Vanishing Theorem now implies that $h^0(K-D)=h^1(K-D)=0$. Serre duality now implies that
$h^2(D)=h^1(D)=0$.
\end{pf}

\begin{lem}
\label{lem.h1D}
For all $n \ge 0$,  $h^1(D_n)=h^2(D_n)=0$.  
\end{lem}
\begin{pf}
The value of $D_n$ for $0 \le n \le 6$ is given explicitly in section \ref{sect.Dn}.  We also note that $D_7 =
D_1+D_6=(4,2,1,2)$. It is routine to check  that conditions  (\ref{eq.ample2}) and  (\ref{eq.ample3}) hold for $D=D_n$ when $n=0,2,3,4,5,6,7$. Hence $h^1(D_n)=h^2(D_n)=0$ when $n=0,2,3,4,5,6,7$. 

We now consider $D_1$ which is the $(-1)$-curve $X=0$. (Since $(D_1-K).D_1 = 0$, $D_1-K$ is not ample so we can't use Kodaira Vanishing as we did for the other small values of $n$.)
 It follows from the exact sequence 
$0 \to \cO_{\BB_3} \to \cO_{\BB_3}(D_1) \to \cO_{D_1}(D_1) \to 0$
that $H^p(\BB_3,\cO_{\BB_3}(D_1)) \cong H^p(\BB_3,\cO_{D_1}(D_1))$ for $p=1,2$. 
However, $D_1 \cong \PP^1$, $\cO_{D_1}(D_1)$ is the
normal sheaf for $D_1\subset \BB_3$, and, since 
$D_1$ can be contracted to a smooth point on the degree 7 del Pezzo
surface, $\cO_{D_1}(D_1) \cong \cO_{D_1}(-1)$. Therefore $H^p(\BB_3, \cO_{D_1}(D_1))
\cong H^p(\PP^1,\cO(-1))$ which is zero for $p=1,2$. It follows that $h^1(D_1)=h^2(D_1)=0$.

Thus  $h^1(D_n)=h^2(D_n)=0$ when $0 \le n \le 7$. We have also shown that
$D_n-K$ is ample when $2 \le n \le 7$.  We now argue by induction. Suppose $n \ge 8$ and $D_{n-6}-K$
is ample. Now $D_n-K=D_{n-6}-K -K$. Since a sum of ample divisors is ample, $D_n-K$ is ample. It follows that $h^1(D_n)=h^2(D_n)=0$.  
\end{pf}

\subsection{The twisted homogeneous coordinate ring $B(\BB_3,\cL,\s)$}
We assume the reader is somewhat familiar with the notion of twisted homogeneous coordinate rings. 
Standard references for that material are \cite{ATV1}, \cite{ATV2}, \cite{AV}, and \cite{AZ}. 

The notion of a $\s$-ample line bundle \cite{AV} plays a key role in the study of twisted homogeneous coordinate rings. Because $\cL_6$  is the anti-canonical bundle and therefore ample, $\cL_1$ is $\s$-ample.
This allows us to use the results of Artin and Van den Bergh in \cite{AV} to conclude that 
the twisted homogeneous coordinate ring
\begin{equation}
\label{defn.B}
B=B(\BB_3,\cL,\s) = \bigoplus_{n=0}^\infty B_n = \bigoplus_{n=0}^\infty H^0(\BB_3,\cL_n)
\end{equation}
is such that 
\begin{equation}
\label{equiv.cats}
\Qcoh \BB_3 \equiv {{\Gr B}\over{\Fdim B}}
\end{equation}
where $\Fdim B$ is the full subcategory of $\Gr B$ consisting of those
graded modules that are the sum of their finite dimensional submodules. 
Artin and Van den Bergh \cite{AV} show that the equivalence (\ref{equiv.cats}) implies that $B$ has 
a host of good properties.

\subsection{}
We will now compute the dimensions $h^0(D_n)$ of the homogeneous  $B_n$ of $B$. 
We will show that $B$ has the same Hilbert series as the non-commutative ring $R$, i.e., the 
same Hilbert series as the weighted polynomial ring with weights 1, 2, and 3. The Hilbert series of $R$
was computed in Proposition \ref{prop.steph}.

As usual  we write
$\chi(D) =h^0(D)-h^1(D)+h^2(D)$.  The Riemann-Roch formula is
$$
\chi(\cO(D))=\chi(\cO)+\hbox{${{1}\over{2}}$} D\cdot (D-K).
$$
 We have $\chi(\cO_{\BB_3})=1$ and $K^2=  6$.

\begin{lem}
\label{lem.h0D}
Suppose $0 \le r \le 5$. Then 
$$
h^0(D_{6m+r})=\begin{cases} (m+1)(3m+r) & \text{if $r \ne 0$,}
\\
3m^2+3m+1 & \text{if $r=0$}
\end{cases}
$$
and 
$$
\sum_{n=0}^\infty h^0(D_n) t^n ={{1}\over{(1-t)(1-t^2)(1-t^3)}} .
$$
In particular, $B$ and $R$ have the same Hilbert series.
\end{lem}
\begin{pf}
Computations for $1 \le r \le 5$ give 
$D_{r}^2=r-2$ and $D_{r}\cdot K= -r $.
Hence
\begin{align*}
\label{T.series}
\chi(D_{6m+r}) &= 1 + \hbox{${{1}\over{2}}$} (D_{r}-mK)\cdot (D_{r}-(m+1) K)
\\
 &= 1 + \hbox{${{1}\over{2}}$} \big(D_r^2-(2m+1)D_r.K +6m(m+1)^2\big)
 \\
 & = (3m+r)(m+1)
\end{align*}
for $m \ge 0$ and $1 \le r \le 5$. When $r=0$, $D_r=0$ so $$\chi(D_{6m}) = 3m^2+3m+1.$$
By Lemma \ref{lem.h1D},  $\chi(D_n)=h^0(D_n)$ for all $n \ge 0$ so
it follows from the formula for $\chi(D_n)$ that
\begin{equation}
\label{eq.chi.diff}
h^0(D_{n+6})-h^0(D_n)=n+6
\end{equation}
 for all $n \ge 0$.

To complete the proof of the lemma, it suffices to show that $h^0(D_n)$ is the coefficient of $t^n$ in the Taylor series expansion 
$$
f(t):={{1}\over{(1-t)(1-t^2)(1-t^3)}} = \sum_{n=0}^\infty  a_n t^n.
$$
 Because
$$
(1-t^6)f(t) = (1-t+t^2)(1-t)^{-2} = 1 + \sum_{n=1}^\infty nt^n,
$$
it follows that
\begin{align*}
(1-t^6)f(t)  & =  a_0+a_1t+\cdots+a_5t^5 + \sum_{n=0}^\infty (a_{n+6}-a_n)t^{n+6} 
\\
& = 1+t+2t+\cdots+5t^5 + \sum_{n=6}^\infty nt^n
\\
& = 1+t+2t+\cdots+5t^5 + \sum_{n=0}^\infty (n+6)t^{n+6}.
\end{align*}
In particular, if $0 \le r \le 5$, $a_r=h^0(D_r)$. 
We now complete the proof by induction. Suppose we have proved that $a_i=h^0(D_i)$ for $i \le n+5$. 
By comparing the expressions in the Taylor series we see that  
$$
a_{n+6}= a_n + (n +6)= h^0(D_n)+n+6= h^0(D_{n+6})
$$
where the last equality is given by (\ref{eq.chi.diff}).
\end{pf}

\subsubsection{Remark.}
It wasn't necessary to compute $\chi(D_n)$ in the previous proof. The proof only used the fact that $\chi(D_{n+6})-\chi(D_n)=n+6$ which can be proved directly as follows:  
\begin{align*}
\chi(D_{n+6})-\chi(D_n)  &=  \hbox{${{1}\over{2}}$} D_{n+6}\cdot(D_{n+6}-K) -  \hbox{${{1}\over{2}}$} D_n\cdot(D_{n}-K) 
\\
& =  \hbox{${{1}\over{2}}$} (D_{n+6}-D_n)\cdot(D_{n+6}+D_n -K)
\\
&=-K\cdot(D_r-(m+1)K)
\\
& = 6(m+1) - K\cdot D_r
\\
&=n+6.
\end{align*}

\subsection{The isomorphism $R \to B(\BB_3,\cL,\s)$}
By definition,  $B_n=H^0(\BB_3,\cL_n)$.
Since Cox's homogeneous coordinate ring, $S=\CC[X,Y,Z,s,t,u]$, is the direct sum of $H^0(\BB_3,\cL)$
as $[\cL]$ ranges over $\Pic \BB_3$ , each $B_n$ is a subspace of $S$. In particular, $B$ itself is a subspace
of $S$, but
\begin{center}
{\it 
  the multiplication in $B$ is not that in $S$.
}
\end{center}

The ring $B$ has the following basis elements in the following degrees:
$$
\begin{array}{ccccc}
\begin{array}{c}
\deg =n \\
1 \\
2 \\
3 \\
4 \\
5 \\
6 \\
\phantom{x}
\end{array}
\begin{array}{c}
\cL_n \\
\cO(1,1,0,1) \\
\cO(1,1,0,0) \\
\cO(2,1,1,1) \\
\cO(2,1,0,1) \\
\cO(3,2,1,1) \\
\cO(3,1,1,1) \\
\phantom{x}
\end{array}
\begin{array}{l} 
\hbox{basis for $B_n$} \\
X \\ 
Xu \\
XYu\\
XYtu\\
XYZtu\\
XYZstu \\
\phantom{x}
\end{array}
\begin{array}{l} 
\\
\\ 
Zt \\
YZt\\
YZt^2\\
YZ^2t^2\\
YZ^2st^2 \\
\phantom{x}
\end{array}
\begin{array}{l} 
\\
\\ 
\\
XZs\\
XZst\\
XZ^2st\\
XZ^2s^2t \\
\phantom{x}
\end{array}
\begin{array}{l} 
\\
\\ 
\\
\\
X^2su\\
X^2Zsu\\
X^2Zs^2u \\
XY^2tu^2 \\
\end{array}
\begin{array}{l} 
\\
\\ 
\\
\\
\\
X^2Yu^2\\
X^2Ysu^2 \\
Y^2Zt^2u.
\end{array}
\end{array}
$$

The multiplication in $B$ is Zhang's twisted multiplication \cite{Ztwist} with respect to the automorphism $\tau$
defined in (\ref{defn.tau}): the product in $B$ of $a \in B_m$ and $b \in B_n$ is
 \begin{equation}
\label{eq.Z.tw.mult}
a *_B b := a\tau^m(b).
 \end{equation}
To make it clear whether a product is being calculated in $B$ or $S$ we will write $x$ for $X$ considered as an element of $B$ and   $y$ for $Zt$ considered as an element of $B$. Therefore, for example, 
\begin{align*}
x^5 & =X\tau(X)\tau^2(X)\tau^3(X)\tau^4(X)\tau^5(X) 
\\
&
= XuYtZ
\\
&= (Zt) Y(uX) 
\\
&
=  Zt\tau^2(X)\tau^3(Zt) 
 \\
 & = yxy
\end{align*}
and 
$$
y^2=Zt\tau^2(Zt)=Zt(sX) = X(sZ)t =  X\tau(zt)\tau^3(X)=xyx.
$$
The following proposition is an immediate consequence of these two calculations.

\begin{prop}
\label{prop.R->B}
Let $R$ be the free algebra $\CC\langle x,y\rangle$ modulo the relations $x^5=yxy$ and $y^2=xyx$.
Then there is a $\CC$-algebra  homomorphism 
$$
\Phi:R=\CC[x,y] \to B(\BB_3,\cL,\s), \qquad x \mapsto X, \; y \mapsto Zt.
$$ 
\end{prop}

\begin{lem}
\label{lem.low.degs}
The homomorphism in Proposition \ref{prop.R->B} is an isomorphism in degrees $\le 6$.\footnote{We will eventually prove that $\Phi$ is an isomorphism in all degrees but the low degree cases need to be handled separately.}
\end{lem}
\begin{pf}
By Proposition \ref{prop.steph}, $R$ has Hilbert series $(1-t)^{-1}(1-t^2)^{-1}(1-t^3)^{-1}$, 
so the dimension of $R_n$ in degrees $1,2,3,4,5,6$ is $1,2,3,4,5,7$. 

The $n^{\th}$ row in the following table gives a basis for $B_n$, $1 \le n \le 6$. One proceeds down each column by multiplying on the right by $x$. There wasn't enough room on a single line for   $B_6$ so we put the last two entries for $B_6$ on a new line.
$$
\begin{array}{cccccc}
\begin{array}{l} 
x=X \\ 
x^2=Xu \\
x^3=XYu\\
x^4=XYtu\\
x^5=XYZtu\\
x^6=XYZstu\\
\phantom{x}
\end{array}
\!\!\!\!
\begin{array}{l} 
\\ 
y=Zt \\
yx=YZt\\
yx^2=YZt^2\\
yx^3=YZ^2t^2\\
yx^4=YZ^2st^2\\
\phantom{x}
\end{array}
\!\!
\begin{array}{l} 
\\ 
\\
xy=XZs\\
y^2=XZst\\
y^2x=XZ^2st\\
y^2x^2=XZ^2s^2t\\
\phantom{x}
\end{array}
\!\!
\begin{array}{l} 
\\ 
\\
\\
x^2y=X^2su\\
xy^2=X^2Zsu\\
xy^2x=X^2Zs^2u\\
 yx^2y=Y^2Zt^2u
\end{array}
\!\!
\begin{array}{l} 
\\ 
\\
\\
\\
x^3y=X^2Yu^2\\
x^2y^2=X^2Ysu^2\\
x^4y=XY^2tu^2.
\end{array}
\end{array}
$$
These calculations involving $x$ and $y$ are made by using the formula (\ref{eq.Z.tw.mult}) in the same way 
it was used to show that $x^5=yxy$.
\end{pf}

\begin{lem}
\label{lem.gend.gl.sects}
$\cL_2$ is generated by its global sections.
\end{lem}
\begin{pf}
A line bundle  on a variety is generated by its global sections if and only if for each point  
 on the variety there is a section of  of the bundle that does not vanish at that point. In this case, $H^0(\BB_3,\cL_2)$ is spanned by $Xu$ and $Zt$. One can see from the diagram (\ref{six.lines}) that the zero locus of $Xu$ does not meet the zero locus of $Zt$, so the common zero locus of $Xu$ and $Zt$ is empty.
\end{pf}

\begin{prop}
\label{prop.B1B2}
As a  $\CC$-algebra, $B$ is generated by $B_1$ and $B_2$.
\end{prop}
\begin{pf}
It follows from the explicit calculations in Lemma \ref{lem.low.degs} that the subalgebra of $B$ generated by
$B_1$ and $B_2$ contains $B_m$ for all $m \le 6$. It therefore 
suffices to prove that the twisted multiplication map  $B_2 \otimes B_n \to B_{n+2}$ is surjective for all 
$n \ge 5$.

By definition, $B_2=H^0(\cL_2)$ and this has dimension two so, by Lemma \ref{lem.gend.gl.sects}, 
there is an exact sequence
$0 \to \cN \to B_2 \otimes \cO_{\BB_3} \to \cL_2 \to 0$ for some line bundle $\cN$. In fact, 
$\cN \cong \cL_2^{-1}$.  

By definition, $\cL_{n+2}=\cL_2 \otimes \cM$ where $\cM \cong \cO(D_{n+2}-D_2)$, and 
the twisted multiplication map $B_2 \otimes B_n \to B_{n+2}$ is  the ordinary multiplication map 
$$
B_2 \otimes H^0(\cM) =  H^0(\cL_2) \otimes H^0(\cM) \to H^0(\cL_2 \otimes \cM).
$$

 Applying $-\otimes \cM$ to the exact sequence $0 \to \cL_2^{-1} \to B_2 \otimes \cO_{\BB_3} \to \cL_2 \to 0$ and taking cohomology gives an exact sequence
 $$
 0 \to H^0(\cL_2^{-1} \otimes \cM) \to B_2 \otimes H^0(\cM) \to H^0(\cL_2 \otimes \cM) \to H^1(\cL_2^{-1} \otimes \cM).
 $$
 Hence, if $h^1(\cL_2^{-1} \otimes \cM)=0$, then $B_2B_n=B_{n+2}$. 
 
 We will now show that $h^1(\cL_2^{-1} \otimes \cM)=0$.
Since
 $
 \cL_2^{-1} \otimes \cM \cong \cO(-D_2 +D_{n+2}-D_2)
 $
 and $n +2 \ge 7$, 
  $$
 \cL_2^{-1} \otimes \cM \cong \cO(-D_{6m+r}-2D_2-K)
 $$
 for suitable integers $m$ and $r$ such that  $6m+r \ge 7$ and $0 \le r \le 5$.
 
By Lemma \ref{lem.Dn}, $D_{6m+r}-2D_2-K =D_r-2D_2-(m+1)K$. By Lemma \ref{lem.van}, 
to show that $h^1(D_r-2D_2-(m+1)K)=0$ it suffices to show that conditions 
(\ref{eq.ample2}) and  (\ref{eq.ample3})  hold for the divisors $D$ in the following table:
 $$
 \begin{array}{lll}
 &&  D:=D_r-2D_2 - (m+1)K \in \Pic\BB_3\\
 r=0   \quad & m \ge 2 \quad & (3m+1,m-1,m+1,m+1) \in \ZZ^4 \\
  r=1 \quad & m \ge 1  \quad& (3m+2,m,m+1,m+2) \\
   r=2  \quad & m \ge 1 \quad & (3m+2,m,m+1,m+1) \\
    r=3  \quad & m \ge 1  \quad& (3m+3,m,m+2,m+2) \\
     r=4  \quad & m \ge 1  \quad & (3m+3,m,m+1,m+2) \\
      r=5  \quad & m \ge 1 \quad & (3m+4,m+1,m+2,m+2).
      \end{array}
$$       
This is a routine task.
\end{pf}

\begin{thm}
\label{thm.isom}
Let $R$ be the free algebra $\CC\langle x,y\rangle$ modulo the relations $x^5=yxy$ and $y^2=xyx$.
The $\CC$-algebra  homomorphism
$$
\Phi: R=\CC[x,y] \to B(\BB_3,\cL,\s), \qquad x \mapsto X, \; y \mapsto Zt,
$$ 
is an isomorphism of graded algebras.
\end{thm}
\begin{pf}
By Lemma \ref{lem.low.degs}, $B_1$ and $B_2$ are in the image of $\Phi$.
By Proposition \ref{prop.B1B2}, $B$ is generated by $B_1$ and $B_2$. 
Hence $\Phi$ is surjective. But $\Phi(R_n) \subset B_n$, 
and $R$ and $B$ have the same Hilbert series, so $\Phi$ is also surjective. 
\end{pf}

Consider $R^{(3)}\supset \CC[x^3,xy,yx]$.
Since $\dim R_6=7=(\dim R_3)^2 -2$ there is a 2-dimensional space of quadratic relations among the elements $x^3$, $xy$, and $yx$. Hence $R^{(3)}$ is not a 3-dimensional Artin-Schelter regular algebra.
 The relations in the degree two component of  $R^{(3)}$ are generated by 
$$
(x^3)^2=(xy)^2=(yx)^2.
$$

 \section{Acknowledgments}
 The author is grateful to the following people: 
 Amer Iqbal for bringing the paper \cite{BP} to his attention; 
 Paul Hacking and  S\'andor Kov\'acs for passing on some standard results about  del Pezzo surfaces
 and algebraic geometry; 
 Darin Stephenson for telling the author that $R$ is an iterated Ore extension and for useful remarks
 about his papers \cite{St1} and \cite{St2}; 
Mark Blunk for pointing out an error in an earlier version of this paper.

\end{document}